\documentclass[11pt]{article}

\begin{document}

\begin{center}
Convolution and Combination Matrices in Non-stationary Filtering \\
Peter  Zizler \\ 
Mount Royal University, Calgary
\end{center}

\begin{abstract}
Time independent convolution yields circulant matrices whose eigenvectors are the Fourier exponentials
with the eigenvalues being the Fourier transform of the mask. The case of time dependent convolution, the non-stationary case, no longer has this property and two matrices are 
then introduced, the cyclic convolution matrix and the cyclic combination matrix. 
In our paper, we prove results on the properties of these matrices. We give results in the context of the non-stationary frequency response in the Fourier domain, where the Fourier 
matrix is a full matrix in general. The techniques used here are attainable at the advanced undergraduate linear algebra settings and can be introduced into a 
relevant linear algebra undergraduate course.
\end{abstract}

\begin{tabular}{ll}
\end{tabular}

\begin{tabular}{ll}
Keywords: &  convolution and combination matrices  \\
~ & non-stationary filtering.
\end{tabular}

\section{Introduction} 

A circulant matrix has the Fourier exponentials as eigenvectors and the Fourier transform of the corresponding mask
are the eigenvalues. This is a foundational result on which the finite Fourier frequency filtering rests. The matrix is created from a given time independent vector, the convolution mask.

In the many application settings, however, the case of time dependent convolution masks arises. In this non-stationary case the Fourier exponentials
are no longer the eigenvectors and the frequency response in general is a full Fourier matrix.
Non-stationary frequency filtering naturally appears in science, in particular in geology, medical imaging, where feature extractions,
trend analysis are implemented. The finite cyclic case for the non-stationary frequency filtering is of uttermost importance as this is the paradigm where relevant algorithms are implemented.
   
We show the non-stationary case involves both the cyclic convolution matrix and the cyclic combination matrix.
The time dependent convolution mask matrix produces the cyclic convolution matrix in the time domain as well as the cyclic
combination matrix in the Fourier domain. We show how these matrices are related and prove some foundational results in their regard.
The stationary case, where the mask is constant and time independent, 
yields the cyclic convolution matrix and the cyclic combination
matrix the same, the classical circulant matrix.

In our paper we give new results on the non-stationary frequency filtering involving the time dependent convolution and combination matrices. 
In particular, we give results when the time dependent mask matrix is of rank one. Results on fast computations are given as well.
We provide a multitude of computational examples to further illustrate the underlying theory. 

The time dependent case, where the mask is no longer stationary, has been studied extensively both in theoretic and applied settings.
In fact, this topic has vast and diverse areas of applications. To this end we refer the reader to  \cite{A}, \cite{M}, \cite{ML}, \cite{MLH}, \cite{R}. 
For foundational results we refer the reader to \cite{D}, \cite{G}.

The time dependent masks are allowed to be complex valued in general, understanding the frequency response matrix will have certain Fourier conjugate symmetries if the masks are real.
Consider the time dependent masks ${\bf c}_{\tau}=\left( c_{0,\tau}, c_{1,\tau}, \ldots, c_{N-1,\tau} \right)^T$ with constant time increments. Without loss of generality we assume unit time 
increments $\{0,1,\ldots,N-1 \}$. We form a (mask) matrix $C$ with columns being the time dependent masks ${\bf c}_{\tau}$, $\tau \in {0,\ldots, N-1}$. For reader's convenience
we provide an example with $N=4$ for the forthcoming preliminary discussions. In particular, the mask matrix $C$ in this case takes the form

$$
C=
\left( \begin{array}{cccc}
c_{00} & c_{01} & c_{02} & c_{03}\\
c_{10} & c_{11} & c_{12} & c_{13}\\
c_{20} & c_{21} & c_{22} & c_{23}\\
c_{30} & c_{31} & c_{32} & c_{33}\\
\end{array}\right)
$$

Based on the matrix $C$ we introduce the cyclic convolution operator

\begin{eqnarray*}
{\bf y}(t) & =&  \sum_{\tau=0}^{N-1} C(t-\tau,\tau) {\bf x}(\tau) \\
\end{eqnarray*}

This operator, in the $N=4$ case, translates to the following matrix multiplication

\begin{eqnarray*}
{\rm conv}(C){\bf x} & = & 
\left( \begin{array}{cccc}
c_{00} & c_{31} & c_{22} & c_{13}\\
c_{10} & c_{01} & c_{32} & c_{23}\\
c_{20} & c_{11} & c_{02} & c_{33}\\
c_{30} & c_{21} & c_{12} & c_{03}\\
\end{array}\right) 
\left( \begin{array}{c}
x_0 \\
x_1 \\
x_2 \\
x_3 \\
\end{array}\right)
\end{eqnarray*}

\noindent where the vector ${\bf x} \in {\bf C}^4$. The general construction of the convolution matrix is given by the following.

$$
{\rm conv}(C)(i,j) = C(i-j,j)
$$

It is not difficult to see the inverse operation is obtained as

$$
C(i,j) = {\rm conv}(C)(i+j,j)
$$

Now the cyclic combination operator is given by 

\begin{eqnarray*}
{\rm comb}(C){\bf x} & =&  \sum_{\tau=0}^{N-1} C(t-\tau,t) {\bf x}(\tau) \\
\end{eqnarray*}

and we observe the case of $N=4$ translates to the following matrix multiplication

\begin{eqnarray*}
{\rm comb}(C){\bf x} & =&  \sum_{\tau=0}^{N-1} C(t-\tau,t) {\bf x}(\tau) \\
& = & 
\left( \begin{array}{cccc}
c_{00} & c_{30} & c_{20} & c_{10}\\
c_{11} & c_{01} & c_{31} & c_{21}\\
c_{22} & c_{12} & c_{02} & c_{32}\\
c_{33} & c_{23} & c_{13} & c_{03}\\
\end{array}\right) 
\left( \begin{array}{c}
x_0 \\
x_1 \\
x_2 \\
x_3 \\
\end{array}\right)
\end{eqnarray*}

The general construction of the combination matrix is given by

$$
{\rm comb}(C)(i,j) = C(i-j,i)
$$

and the inverse operation is obtained as

$$
C(i,j) = {\rm comb}(C)(j,j-i)
$$

The following result connects the convolution and combination matrix directly without the immediate reference to the mask matrix $C$.

\noindent {\bf Proposition 1.1}
We have 

$$
{\rm comb}(C)(i,j)={\rm conv}(C)(2i-j,i)
$$

$$
{\rm conv}(C)(i,j)={\rm comb}(C)(j,2j-i)
$$

{\bf Proof} Note

$$
{\rm conv}^{-1} \left( {\rm comb}(i,j) \right) = {\rm conv}^{-1} \left(i-j,i \right) = (i-j+i,i) = (2i-j,i)
$$

and

$$
{\rm comb}^{-1} \left({\rm conv}(i,j) \right) = {\rm comb}^{-1} \left(i-j,j \right) = (j,j-(i-j) = (j,2j-i)
$$

We observe that in the stationary case we have ${\bf c}_{\tau} = {\bf c}=(c_0, c_1, \ldots, c_{N-1})^T$. In particular, we have one time independent mask ${\bf c}$.
Continuing with our $N=4$ example setting we note the stationary (mask) matrix 

$$
C=
\left( \begin{array}{cccc}
c_0 & c_{0} & c_{0} & c_{0}\\
c_1 & c_{1} & c_{1} & c_{1}\\
c_2 & c_{2} & c_{2} & c_{2}\\
c_3 & c_{3} & c_{3} & c_{3}\\
\end{array}\right)
=
\left( \begin{array}{c}
c_0 \\
c_1 \\
c_2 \\
c_3 \\
\end{array}\right)
\left( \begin{array}{cccc}
1 & 1 &  1 & 1 
\end{array}\right)
$$

\noindent with 

$$
{\rm conv}(C)={\rm comb}(C) = \left( \begin{array}{cccc}
c_0 & c_{3} & c_{2} & c_{1}\\
c_1 & c_{0} & c_{3} & c_{2}\\
c_2 & c_{1} & c_{0} & c_{3}\\
c_3 & c_{2} & c_{1} & c_{0}\\
\end{array}\right)
$$

It is straightforward to see that in the stationary case in general we have ${\rm conv}(C)={\rm comb}(C)$. This matrix is referred to as the circulant matrix. 

The discrete Fourier transform of the vector ${\bf x}$ is given by

$$
{\bf X}(j)={\rm DFT}({\bf x})(j) = \sum^{N-1}_{k=0} {\bf x}(k) e^{\frac{-2 \pi i j k}{N}} 
$$

\noindent and the discrete inverse Fourier transform of the vector ${\bf X}=(X_0, X_1, \ldots, X_{N-1})^T$ is given by

$$
{\bf x}(k) =  \frac{1}{N}\sum^{N-1}_{j=0} {\bf X}(j) e^{\frac{2 \pi i j k}{N}} 
$$

The two dimensional discrete Fourier transform of the (mask) matrix $C$ is defined as

$$
G(p,q) = {\rm DFT2}(C)(p,q) =  \sum^{N-1}_{m=0} \sum^{N-1}_{n=0} C(m,n) e^{\frac{-2 \pi i p m}{N}} e^{\frac{-2 \pi i q n}{N}}   
$$

with the discrete inverse two dimensional Fourier transform defined as 

$$
C(m,n) = {\rm IDFT2}(G)(m,n) = \frac{1}{N^2}\sum^{N-1}_{p=0} \sum^{N-1}_{q=0} G(p,q) e^{\frac{2 \pi i p m}{N}} e^{\frac{2 \pi i q n}{N}} 
$$

We define a $N \times N$ unitary (Fourier) matrix 

$$
V=\frac{1}{\sqrt{N}}\left[ v_{k,j}\right]
$$

where $v_{k,j} = e^{\frac{2 \pi i j k}{N}}$. We observe that

$$
{\rm DFT}({\bf x}) = \sqrt{N} \,  V^{*}{\bf x} \mbox{ ; } {\rm IDFT}({\bf x})= \frac{1}{\sqrt{N}} \, V{\bf x}
$$

$$
{\rm DFT}\left( {\rm IDFT}({\bf x}) \right) = {\bf x} \mbox{ ; } {\rm IDFT}\left( {\rm DFT}({\bf x}) \right) = {\bf x}
$$

$$
VCV=N \, {\rm IDFT2}(C) \mbox{ ; } V^{*}CV^{*}=\frac{1}{N} \, {\rm DFT2}(C) 
$$

\section{\bf Preliminaries}

We define a linear operator $J$ by the following rule: $x_0 \mapsto x_0$ and $x_k \mapsto x_{N-k}$ \mbox{ ; } with ${\bf x} = (x_0,x_1, \ldots, x_{N-k})^T$. 
In the $N=4$ example setting we have

$$
J = \left( \begin{array}{cccc}
1 & 0 & 0 & 0\\
0 & 0 & 0 & 1\\
0 & 0 & 1 & 0\\
0 & 1 & 0 & 0\\
\end{array}\right)
$$

In general we have the following

$$
\left( V^{*}\right)^2{\bf x} = \frac{1}{N}{\rm DFT}\left( {\rm DFT}({\bf x}) \right)  =J{\bf x}
$$

$$
V^2{\bf x} = N \, {\rm IDFT}\left( {\rm IDFT}({\bf x}) \right)  =J{\bf x}
$$

We now define $F=V^*C^TV^*$. The following result is taken from \cite{Z}. For reader's convenience we include a proof in the context of the above setting.

\noindent {\bf Lemma 2.1} We have

$$
V^{*}{\rm conv}(C) V = {\rm comb}(F) \mbox{ ; } V^{*}{\rm comb}(C) V = {\rm conv}(F) 
$$

{\bf Proof} Let ${\bf x} \in {\bf C}^n$ and consider

\begin{eqnarray*}
{\rm DFT}\left({\rm comb}(C){\bf x}\right) & = & \sum_t e^{\frac{-2 \pi i j t}{n}} \sum_{\tau} C_{\rm comb}(t-\tau,t) {\bf x}(\tau)\\
& = & \frac{1}{n} \sum_t e^{\frac{-2 \pi i j t}{n}} \sum_{\tau} C_{\rm comb}(t-\tau,t) \sum_s {\bf X}(s) e^{\frac{2 \pi i \tau s}{n}} \\
& = & \frac{1}{n} \sum_s \sum_t \sum_{\tau} C(t -\tau,t ) e^{\frac{2 \pi i \tau s}{n}} e^{\frac{-2 \pi i j t}{n}} {\bf X}(s)
\end{eqnarray*}

\noindent now let $u = t-\tau$ and we continue

\begin{eqnarray*}
{\rm DFT}\left({\rm comb}(C){\bf x}\right) & = & \frac{1}{n} \sum_s \left\{  \sum_t \sum_u C(u,t ) e^{\frac{-2 \pi i u s}{n}} e^{\frac{-2 \pi i (j-s) t}{n}}   \right\} {\bf X}(s) \\
& = & \sum_s {\rm conv}(F) (j-s,s){\bf X}(s) \\
& = & {\rm conv}(F){\bf X}.
\end{eqnarray*}

Now for ${\rm conv}(C)$ the case is very similar

\begin{eqnarray*}
{\rm DFT}\left({\rm conv}(C){\bf x}\right)  & = & \sum_t e^{\frac{-2 \pi i j t}{n}} \sum_{\tau} C_{\rm conv}(t-\tau,\tau) {\bf x}(\tau)\\
& = & \frac{1}{n} \sum_t e^{\frac{-2 \pi i j t}{n}} \sum_{\tau} C_{\rm conv}(t-\tau,\tau) \sum_s {\bf X}(s) e^{\frac{2 \pi i \tau s}{n}} \\
& = & \frac{1}{n} \sum_s \sum_t \sum_{\tau} C(t -\tau,\tau ) e^{\frac{2 \pi i \tau s}{n}} e^{\frac{-2 \pi i j t}{n}} {\bf X}(s)
\end{eqnarray*}

\noindent now let $u = t-\tau$ and we continue

\begin{eqnarray*}
{\rm DFT}\left({\rm conv}(C){\bf x}\right)  & = & \frac{1}{n} \sum_s \left\{  \sum_{\tau} \sum_u C(u,\tau ) e^{\frac{-2 \pi (j-s)\tau}{n}} e^{\frac{-2 \pi i j u}{n}}   \right\} {\bf X}(s) \\
& = & \sum_s {\rm comb}(F) (j-s,j){\bf X}(s) \\
& = & {\rm comb}(F){\bf X}.
\end{eqnarray*}

\section{Results}

Suppose the (mask) matrix $C$ is a Kronecker product $C={\bf c}{\bf d}^{*}$ with ${\bf c},{\bf d}$ complex vectors. We provide some results for fast computations involving this specialized matrix.
Suppose the two vectors ${\bf x}, {\bf y}$ are given. By the notation ${\bf x}.*{\bf y}$ we understand the point wise multiplication ${\bf x}.*{\bf y}(i)={\bf x}(i){\bf y}(i)$.

\noindent {\bf Lemma 3.1} Suppose the matrix $C={\bf c}{\bf d}^{*}$. Then we obtain

$$
{\rm comb}(C){\bf x} = {\rm conj}({\bf d}).*{\rm IDFT}\left( {\rm DFT}({\bf c}).* {\rm DFT}({\bf x}) \right)
$$

{\bf Proof} We have 

\begin{eqnarray*}
{\rm comb}(C){\bf x}(t) & = & \sum_{\tau} C(t-\tau,t){\bf x}(\tau) \\
& = &  \sum_{\tau} {\bf c}(t-\tau){\bf d}^{*}(t){\bf x}(\tau) \\
& = &  {\bf d}^{*}(t) \sum_{\tau} {\bf c}(t-\tau){\bf x}(\tau) \\
& = & {\rm conj}({\bf d}).*{\rm IDFT}\left( {\rm DFT}({\bf c}).* {\rm DFT}({\bf x}) \right) \\
\end{eqnarray*}

\noindent {\bf Lemma 3.2} Suppose the matrix $C={\bf c}{\bf d}^{*}$. Then we obtain

$$
{\rm conv}(C){\bf x} = {\rm IDFT}\left(   {\rm DFT}\left({\rm conj}({\bf d}).*{\bf x}  \right).* {\rm DFT}({\bf c})        \right)
$$

{\bf Proof} We have

\begin{eqnarray*}
{\rm DFT}\left({\rm conv}(C){\bf x}  \right) (j) & = & \sum_t \sum_{\tau} C(t-\tau,\tau){\bf x}(\tau) e^{\frac{-2 \pi j t}{N}  }\\
& = & \sum_t \sum_{\tau} {\bf c}(t-\tau){\bf d}^{*}(\tau){\bf x}(\tau) e^{\frac{-2 \pi j t}{N}  }\\
& = & \sum_{\tau} {\bf d}^{*}(\tau){\bf x}(\tau) \sum_t {\bf c}(t-\tau) e^{\frac{-2 \pi j t}{N}  }\\
& = & \sum_{\tau} {\bf d}^{*}(\tau){\bf x}(\tau) e^{\frac{-2 \pi j \tau}{N}  } \sum_s {\bf c}(s) e^{\frac{-2 \pi j s}{N}  }    \hspace{2cm}  s = t-\tau \\
& = & {\rm DFT}({\bf \lambda})(j){\rm DFT}({\bf c})(j)   \hspace{2cm} {\bf \lambda} = {\rm conj}({\rm d}).*{\bf x}
\end{eqnarray*}

\noindent {\bf Lemma 3.3} Suppose the matrix $C={\bf c}{\bf d}^{*}$. Then we obtain

$$
F = \frac{1}{N} {\rm DFT}({\rm conj}(d)).*{\rm DFT}({\bf c})^T
$$

{\bf Proof} We have

\begin{eqnarray*}
V^* C^T V^* & = & V^* ({\bf c}{\bf d}^*)^T V^* \\
& = & V^* {\rm conj}({\bf d}){\bf c}^T V^* \\
& = & V^* {\rm conj}({\bf d})\left( V {\rm conj}({\bf c})\right)^*\\
& = & \frac{1}{\sqrt{N}} \, {\rm DFT}({\rm conj}({\bf d})) \, \sqrt{N} \,{\rm IDFT}({\rm conj}({\bf c}))^* \\
& = & {\rm DFT}({\rm conj}({\bf d})) \, {\rm IDFT}({\rm conj}({\bf c}))^* \\
& = & \frac{1}{N} \, {\rm DFT}({\rm conj}({\bf d})) \, {\rm DFT}({\bf c})^T \\
\end{eqnarray*}

since

$$
{\rm conj}\left(   {\rm IDFT} \left( {\rm conj}({\bf c}) \right)   \right) = \frac{1}{N} {\rm DFT}({\rm c})
$$

\noindent {\bf Corollary 3.4} Suppose the matrix $C={\bf c}{\bf d}^{*}$. Then we obtain

$$
V^* {\rm conv}(C) V (i,j)= {\rm comb}(F)(i,j) = \frac{1}{N} {\rm DFT}({\rm conj}({\bf d}))(i-j).* {\rm DFT}({\bf c})^T(i)
$$

$$
V^* {\rm comb}(C) V (i,j)= {\rm conv}(F)(i,j) = \frac{1}{N} {\rm DFT}({\rm conj}({\bf d}))(i-j).* {\rm DFT}({\bf c})^T(j)
$$

\section{\bf  Frequency Filtering}

The following observation is useful in our further analysis. 
Let $A$ be a $N \times N$ matrix. Suppose $\{{\bf v}_k\}_{k=0}^{N-1}$ is a given orthonormal set in ${\bf C^N}$. 
Then $A=\sum_{k} A{\bf v}_k {\bf v}^*_k$. Note that this expansion is unique, in a sense that if  $A=\sum_{k} {\bf u}_k {\bf v}^*_k = \sum_{k} {\bf w}_k {\bf v}^*_k$ then 
${\bf w}_k = {\bf u}_k = A{\bf v}_k$. Also

$$
\left< A{\bf v}_k{\bf v}^*_k, A{\bf v}_l{\bf v}^*_l\right>_F = ||A{\bf v}_k||^2 \delta_{k,l} 
$$

where $\left< A,B \right>_F={\rm tr}(B^*A)$ denotes the Frobenious inner product.

Let $C$ be an $N \times N$ real matrix of columns being the time dependent masks. Given the matrix $C$ we obtain $F=\frac{1}{N}{\rm DFT2}(C^T)$ and conversely given the 
matrix $F$ we obtain $C=\left( {\rm IDFT2}(N \, F) \right)^T$. 
It is not difficult to observe that the matrix $C$ is real if and only if

$$
F(i,j)= {\rm conj}(F(i_1,j_1))
$$

whenever $i+i_1 \equiv 0$ mod $N$ and $j+j_1 \equiv 0$ mod $N$, referring to such as the Fourier conjugate symmetry.

\noindent {\bf Proposition 4.1} Suppose the (mask) matrix $C$ is a real valued and consider the corresponding matrix $F=\frac{1}{N}{\rm DFT2}(C^T)$. The following are equivalent.

\begin{itemize}

\item The matrix $F$ has the Fourier conjugate symmetry 

\item The matrix ${\rm comb}(F)$ has the Fourier conjugate symmetry 

\item The matrix ${\rm conv}(F)$ has the Fourier conjugate symmetry

\end{itemize}

{\bf Proof} Suppose $(i_1,j_1)$ is the additive inverse of $(i,j)$ in ${\bf Z}_{N \times N}$. Let $m,n, k, l$ be arbitrary whole numbers. Then 
$(mi_1+nj_1, k i_1+lj_1)$ is the additive inverse of $(mi+nj, ki+lj)$.

Suppose the matrix $F$ is given. Let ${\bf f}_k$ denote the $k$ the row of the matrix $F$. Another worthwhile observation is the following. The matrix $F$ displays the Fourier conjugate symmetry if and only if for $0 \leq k \leq {\rm floor}(\frac{N}{2})$ we have ${\bf f}_{N-k} = {\rm conj}(J({\bf f}_k))$. 

Define $F_0$ as a zero $N \times N$ matrix with the exception of the first row equal to the first row of $F$.
For $0 \leq k \leq {\rm floor}(\frac{N}{2})$ define $F_k$ as a zero matrix with the exception of the $k$ row being the $k$ row of $F$ and the $N-k$ row being the $N-k$ row of $F$.

Suppose the matrix $F$ is given and set $k=i-j$. Let $0 \leq k \leq {\rm floor}(\frac{N}{2})$ then we have ${\rm comb}(F_k)(i,i-k)={\bf f}_k(i)$ 
and ${\rm comb}(F_k)(i,i+k)={\rm conj}(J{\bf f}_k)(i)$. The vector ${\bf e}_k$ denotes the vector of all entries zero except the entry one in the $k$th position. 

\noindent {\bf Lemma 4.2} Let ${\bf exp}_j(k) = e^\frac{2 \pi i j k }{N}$. Suppose the matrix $F$ is given. We set $C=\left( {\rm IDFT2}(N \, F) \right)^T$. Define the matrix $F_k$
for $0 \leq k \leq {\rm floor}(\frac{N}{2})$. Set $C_k=\left( {\rm IDFT2}(N \, F_k) \right)^T$. We have

\begin{itemize}

\item 

\begin{eqnarray*}
C_0 & = & {\rm IDFT}({\bf f}_0) {\bf exp}^T_0  \\
\end{eqnarray*}

\begin{eqnarray*}
C_0 & = & \frac{1}{N} C \left( \begin{array}{c}
1 \\
1\\
\vdots \\
1\\
\end{array}\right)
\left( \begin{array}{cccc}
1 & 1 & \dots & 1 \\
\end{array}\right)
\end{eqnarray*}

\item For $0 \leq k \leq {\rm floor}(\frac{N}{2})$

\begin{eqnarray*}
C_k & = & 2 \, {\rm Re}\left( {\rm IDFT}({\bf f}_k) {\bf exp}^T_j \right) \\
\end{eqnarray*}

\begin{eqnarray*}
C_k & = & \frac{1}{N} \left( C{\bf exp}_j{\bf exp}_j^* + C{\bf exp}_{N-j}{\bf exp}_{N-j}^*  \right) \\
& = & \frac{1}{N}  C \left( \begin{array}{c}
1 \\
e^{\frac{2 \pi i j (1)}{N}}\\
\vdots \\
e^{\frac{2 \pi i j (N-1)}{N}}\\
\end{array}\right)
\left( \begin{array}{cccc}
1 & e^{\frac{-2 \pi i j (1)}{N}} & \dots & e^{\frac{-2 \pi i j (N-1)}{N}} \\
\end{array}\right)  \\
 & + & \frac{1}{N}  
C \left( \begin{array}{c}
1 \\
e^{\frac{-2 \pi i j (1)}{N}}\\
\vdots \\
e^{\frac{-2 \pi i j (N-1)}{N}}\\
\end{array}\right)
\left( \begin{array}{cccc}
1 & e^{\frac{2 \pi i j (1)}{N}} & \dots & e^{\frac{2 \pi i j (N-1)}{N}} \\
\end{array}\right)
\end{eqnarray*}

\item If $N=2n$

\begin{eqnarray*}
C_n & = & {\rm IDFT}({\bf f}_n){\bf exp}^T_n \\
\end{eqnarray*}

\begin{eqnarray*}
C_n & = & \frac{1}{N} C \left( \begin{array}{c}
1 \\
-1\\
\vdots \\
-1\\
\end{array}\right)
\left( \begin{array}{cccc}
1 & -1 & \dots & -1 \\
\end{array}\right)
\end{eqnarray*}

\end{itemize}

In particular the frequency response of ${\rm conv}(C_k)$ is ${\rm comb}(F_k)$.

{\bf Proof} Note

\begin{eqnarray*}
C_0 & = & \left( N \, {\rm IDFT2}\left( {\bf e}_0 {\bf f}^T_0\right) \right)^T \\
& = & \left( V {\bf e}_0 {\bf f}^T_0 V \right)^T \\
& = & \left( {\rm IDFT}({\bf e}_0 ) {\rm DFT}(({\rm conj}({\bf f}_0))^*)   \right)^T \\
& = & {\rm IDFT}({\bf f}_0) {\bf exp}^T_0 \\
& = & {\rm IDFT}({\bf f}_0)\left(1 , 1 , 1,  \ldots, 1 \right)
\end{eqnarray*}

Case $N=2n$

\begin{eqnarray*}
C_n & = & \left( N \, {\rm IDFT2}\left( {\bf e}_n {\bf f}^T_n\right) \right)^T \\
& = & \left( V {\bf e}_n {\bf f}^T_n V \right)^T \\
& = & \left( {\rm IDFT}({\bf e}_n ) {\rm DFT}(({\rm conj}({\bf f}_n))^*)   \right)^T \\
& = & {\rm IDFT}({\bf f}_0) {\bf exp}^T_n \\
& = & {\rm IDFT}({\bf f}_n)\left(1 , -1 , 1,  \ldots, -1 \right)
\end{eqnarray*}

For $0 \leq k \leq {\rm floor}(\frac{N}{2})$

\begin{eqnarray*}
C_k & = & \left( N \, {\rm IDFT2}\left( {\bf e}_k {\bf f}^T_k + {\bf e}_{N-k} (J{\bf f}_k)^* \right) \right)^T \\
& = & \left( V {\bf e}_k {\bf f}^T_k V  + V {\bf e}_{N-k} (J{\bf f}_k)^* V  \right)^T \\
& = & \left( {\rm IDFT}({\bf e}_k){\rm DFT}({\rm conj}({\bf f}_k))^* + {\rm IDFT}({\bf e}_{N-k}){\rm DFT}(J{\bf f}_k)^* \right)^T \\
& = & {\rm IDFT}({\bf f}_k){\bf exp}^T_k + {\rm conj}\left( {\rm IDFT}({\bf f}_k) \right) {\bf exp}^T_{N-k} \\ 
& = & 2{\rm Re} \left( {\rm IDFT}({\bf f}_k) {\bf exp}^T_k \right) \\
\end{eqnarray*}

The remaining results follow using the observation.

For reader's convenience we now provide several computational examples in order to illustrate the above results that might be technical in nature.  

\noindent {\bf Example 4.3}

$$
F= \left( \begin{array}{cccc}
 1  &   2-i  &   1   &  2+i \\
-1+2i  &   1-4i  &  -1   &  3+i \\
 3  &   1-i  &   2  &   1+i \\ 
 -1-2i  &   3-i &    -1  &   1+4i \\
\end{array}\right)
\mbox{ ; }
C =\left( \begin{array}{cccc}
4.25  &  0.25  &  2.25  & -0.75 \\
    3.75  & -0.25  & -1.25  & -0.25 \\
   -2.75 &  -3.75  &  3.25  &  1.25 \\
   -3.25   &-2.25  &  1.75  &  1.75 \\
\end{array}\right)
$$

$$
{\rm conv}(C) =\left( \begin{array}{cccc}
 4.25  & -2.25  &  3.25  & -0.25 \\
    3.75 &   0.25 &   1.75 &   1.25 \\
   -2.75 &  -0.25 &   2.25 &   1.75 \\
   -3.25 &  -3.75 &  -1.25 &  -0.75 \\
\end{array}\right)
\mbox{ ; }
F_0= \left( \begin{array}{cccc}
 1  &   2-i  &   1   &  2+i \\
0  &   0  &  0   &  0 \\
 0  &   0  &   0  &   0 \\ 
 0  &   0 &    0  &   0 \\
\end{array}\right)
$$

$$
{\rm comb}(F_0)= \left( \begin{array}{cccc}
 1  &   2-i  &   1   &  2+i \\
0  &   2-i  &  0   &  0 \\
 0  &   0  &   1  &   0 \\ 
 0  &   0 &    0  & 2-i   \\
\end{array}\right)
\mbox{ ; }
F_1= \left( \begin{array}{cccc}
0  &   0  &  0   &  0 \\
-1+2i  &   1-4i  &  -1   &  3+i \\
 0  &   0  &   0  &   0 \\ 
 -1-2i  &   3-i &    -1  &   1+4i \\   \\
\end{array}\right)
$$

$$
{\rm comb}(F_1)= \left( \begin{array}{cccc}
0 & -1 - 2i &  0 & -1 + 2i \\
   1 - 4i &  0 &  3- i &  0 \\
   0 & -1 &  0 & -1  \\
   1 + 4i &  0 &  3 + i &  0 \\
\end{array}\right)
\mbox{ ; }
F_2= \left( \begin{array}{cccc}
0  &   0  &  0   &  0 \\
 0  &   0  &   0  &   0 \\ 
3  &   1-i  &   2  &   1+i \\ 
 0  &   0 &    0  &    \\
\end{array}\right)
$$

$$
{\rm comb}(F_2)= \left( \begin{array}{cccc}
0 &  0  & 3  & 0 \\
   0 &  0 &  0 &  1 - i \\
   2 &  0 &  0 &  0 \\
   0 &  1 + i &  0 &  0 \\
\end{array}\right)
\mbox{ ; }
C_0= \left( \begin{array}{cccc}
  1.5  &  1.5  &  1.5  &  1.5 \\
    0.5 &   0.5  &  0.5 &   0.5 \\
   -0.5 &  -0.5  & -0.5  & -0.5 \\
   -0.5 &  -0.5 &  -0.5 &  -0.5 \\
\end{array}\right)
$$

$$
{\rm conv}(C_0)= \left( \begin{array}{cccc}
   1.5  & -0.5  & -0.5 &   0.5 \\
    0.5 &   1.5  & -0.5 &  -0.5 \\
   -0.5 &   0.5  &  1.5 &  -0.5 \\
   -0.5  & -0.5  &  0.5 &   1.5 \\
\end{array}\right)
\mbox{ ; }
C_1= \left( \begin{array}{cccc}
    1.0 &   0.5  & -1.0  & -0.5 \\
    2.5 &        0.0  & -2.5  & 0.0 \\
   -3.0 &  -2.5  &  3.0  &  2.5 \\
   -2.5 &  -2.0  &  2.5  &  2.0 \\

\end{array}\right)
$$

$$
{\rm conv}(C_1)= \left( \begin{array}{cccc}
      1.0  & -2.0 &   3.0  &       0.0 \\
    2.5  &  0.5  &  2.5  &  2.5 \\
   -3.0   &      0.0 &  -1.0  &  2.0 \\
   -2.5  & -2.5 &  -2.5  & -0.5 \\
\end{array}\right)
\mbox{ ; }
C_2= \left( \begin{array}{cccc}
     1.75 &  -1.75  &  1.75 &  -1.75 \\
    0.75 &  -0.75  &  0.75  & -0.75 \\
    0.75 &  -0.75 &   0.75  & -0.75 \\
   -0.25  &  0.25 &  -0.25  &  0.25 \\
\end{array}\right)
$$

$$
{\rm conv}(C_2)= \left( \begin{array}{cccc}
    1.75  &  0.25  &  0.75 &  -0.75 \\
    0.75  & -1.75 &  -0.25 &  -0.75 \\
    0.75  & -0.75 &   1.75 &   0.25 \\
   -0.25  & -0.75 &   0.75 &  -1.75 \\
\end{array}\right)
$$

\noindent {\bf Example 4.4}

Set 

$$
{\bf exp}_1 =\left( \begin{array}{c}
   1  \\
   \frac{1}{2} + \frac{\sqrt{3}}{2}i \\
  -\frac{1}{2} + \frac{\sqrt{3}}{2}i \\
  -1  \\
  -\frac{1}{2} - \frac{\sqrt{3}}{2}i \\
   \frac{1}{2} - \frac{\sqrt{3}}{2}i \\
\end{array}\right)
$$

and define

$$
C_1 =2 \, {\rm real}\left( {\rm IDFT}\left((1 , 1 , 1 , 1 , 1 , 1)^T\right) {\bf exp}_1^T \right)
$$

we obtain 

$$
C_1 = \left( \begin{array}{cccccc}
    2 &   1 &  -1 & -2  & -1 &   1 \\
         0     &    0     &    0     &    0     &    0     &    0 \\
         0     &    0     &    0     &    0     &    0     &    0 \\
         0     &    0     &    0     &    0     &    0    &     0 \\
         0     &    0     &    0     &    0     &    0   &      0 \\
         0     &    0     &    0     &    0     &    0   &      0 \\
\end{array}\right)
\mbox{ ; }
{\rm conv}(C_1) = \left( \begin{array}{cccccc}
         2     &    0     &    0     &    0     &    0     &   0 \\
         0     &    1     &    0     &    0     &    0     &   0 \\
         0     &    0     &   -1     &    0     &    0     &   0 \\
         0     &    0     &    0     &   -2     &    0     &   0 \\
         0     &    0     &    0     &    0     &   -1     &   0 \\
         0     &    0     &    0     &    0     &    0     &   1 \\
\end{array}\right)
$$

with

$$
F_1 =\left( \begin{array}{cccccc}
  0 & 0 & 0 & 0 & 0  & 0  \\
	1 & 1 & 1 & 1 & 1  & 1  \\
	0 & 0 & 0 & 0 & 0  & 0  \\
	0 & 0 & 0 & 0 & 0  & 0  \\
	0 & 0 & 0 & 0 & 0  & 0  \\
	1 & 1 & 1 & 1 & 1  & 1  \\
\end{array}\right)
\mbox{ ; }
{\rm comb}(F_1) =\left( \begin{array}{cccccc}
  0 & 1 & 0 & 0 & 0  & 1  \\
	1 & 0 & 1 & 0 & 0  & 0  \\
	0 & 1 & 0 & 1 & 0  & 0  \\
	0 & 0 & 1 & 0 & 1  & 0  \\
	0 & 0 & 0 & 1 & 0  & 1  \\
	1 & 0 & 0 & 0 & 1  & 0  \\
\end{array}\right)
$$

\noindent {\bf Lemma 4.5} Let $F$ be a zero matrix except the $k$th row being ${\bf f}_k$ and the $N-k$ row being ${\rm conj}(J{\bf f}_k)$. For any vector ${\bf x} \in {\bf R}^N$ we have

\begin{itemize}

\item If $k=0$

$$
{\rm conv}(C){\bf x} =  {\rm IDFT} \left( {\rm DFT}({\bf x}). * {\bf f}_0 \right)
$$

\item If $N=2n$ and $k=n$

$$
{\rm conv}(C){\bf x} =  {\rm IDFT} \left( {\rm DFT}({\bf x})_n. * {\bf f}_n \right)
$$

where ${\rm DFT}({\bf x})_n(j) = {\rm DFT}({\bf x})(j-n)$

\item otherwise

$$
{\rm conv}(C){\bf x} = 2 \, {\rm Re} \, \left({\rm IDFT} \left( {\rm DFT}({\bf x})_k . * {\bf f}_k \right) \right)
$$

where ${\rm DFT}({\bf x})_k(j) = {\rm DFT}({\bf x})(j-k)$

\end{itemize}

{\bf Proof}

Set ${\bf c}={\rm IDFT}({\bf f}_k)$ and ${\bf d}={\bf exp}_{N-k}$. We have

\begin{eqnarray*}
C{\bf x} & = & {\bf c}{\bf d}^*{\bf x}+ {\rm conj}\left({\bf c}{\bf d}^* \right){\bf x}  
\end{eqnarray*}

If $k=0$ we have

\begin{eqnarray*}
{\rm conv}(C){\bf x} & = & {\rm IDFT}\left( {\rm DFT}({\bf x}).* {\rm DFT}\left( {\rm IDFT}({\bf f}_0)\right) \right) \\
& = &  {\rm IDFT}\left( {\rm DFT}({\bf x}).* {\bf f}_0 \right) \\
\end{eqnarray*}

If $k=n$, where $N=2n$, we have

\begin{eqnarray*}
{\rm conv}(C){\bf x} & = & {\rm IDFT}\left( {\rm DFT}({\bf exp}_n .* {\bf x}).* {\bf f}_k \right) \\
\end{eqnarray*}

otherwise

\begin{eqnarray*}
{\rm conv}(C){\bf x} & = & {\rm IDFT}\left( {\rm DFT}({\bf exp}_k .* {\bf x}).* {\rm DFT}\left( {\rm IDFT}({\bf f}_k)\right) \right) \\
& + &  {\rm conj}\left({\rm IDFT}\left( {\rm DFT}({\bf exp}_k .* {\bf x}).* {\rm DFT}\left( {\rm IDFT}({\bf f}_k)\right) \right) \right)\\
& = &  2 \, {\rm Re} \, \left({\rm IDFT} \left( {\rm DFT}({\bf x})_k . * {\bf f}_k \right) \right) \\
\end{eqnarray*}

\noindent {\bf Example 4.6} We set

$$
{\bf f}_0=(1, 2+i, 3-i, i, 2, -i, 3+i, 2-i)^T \mbox{ ; } {\bf f}_1=(2+3i, 1, 2, -i, i, 3, -1-i, -1)^T 
$$

$$
{\bf f}_2=(-i, i, 1, 4, -i, 2+i, 2+i ,-1 )^T
$$

\noindent understanding 

$$
{\rm comb}(F) =\left( \begin{array}{cccccccc}
    1  &  2-3i &  i  & 0  &  0 &  0  &  -i  &  2+3i \\
    1  &  2+i &  -1  & -1  &  0 &  0  &  0  &  i \\
		1  &  2 &  3-i  & -1+i  &  2-i &  0  &  0  &  0 \\
		0  &  4 &  -i  & i  &  3 &  2-i  &  0  &  0 \\
		0  &  0 &  -i  & i  &  2 &  -i  &  i  &  0 \\
		0  &  0 &  0  & 2+i  &  3 &  -i  &  i  &  4 \\
		1  & 0 &  0  & 0 &  2+i &  -1-i  &  3+i  &  2 \\
		 1  &  -i &  0  & 0  &  0 &  -1  &  -1  &  2-i \\
\end{array}\right)
$$

Let 

$$
{\bf x}=(1 , -2 , 3 , 1 , 1 , 0 , -2 , 1)^T
$$

\noindent and obtain

\begin{eqnarray*}
{\rm conv}(C){\bf x} & = & {\rm IDFT}\left( {\rm DFT}({\bf x}).* {\bf f}_0 \right) \\
& + &  2 \, {\rm Re} \, \left({\rm IDFT} \left( {\rm DFT}({\bf x})_1 . * {\bf f}_1 \right) \right) \\
& + &  2 \, {\rm Re} \, \left({\rm IDFT} \left( {\rm DFT}({\bf x})_2 . * {\bf f}_2 \right) \right) \\
& = & 
(  1.4608  , -0.5895  ,  1.8750  ,  2.3750  ,  4.2892 ,  -5.6605  , -3.1250  ,  2.3750)^T \\
& + & 
( -3.9142  ,  0.1893  , -1.5858 ,  -2.0074  , -11.4142  ,  2.3107  , -3.0858 ,  -7.6642)^T \\
& + & 
(-1.3107  ,  5.9372  , -7.4393  ,  2.9372  ,  5.8107  , -8.4372  ,  2.9393  , -8.4372)^T \\
& = &
( -3.7641   , 5.5371   , -7.1501   , 3.3048  , -1.3143 , -11.7871  , -3.2714 , -13.7264 )^T\\
\end{eqnarray*}

Note ${\rm conv}(C)$ equals

$$
\left( \begin{array}{cccccccc}
    5.1250  &  0.5695 &  -0.1250  & -1.9660  &  0.1250 &  -3.3195  &  3.3750  &  1.2160 \\
   -0.9660  &  2.0821 &   1.7160 &  -1.2286 &   1.0089  &  2.0821  & -2.6731  &  0.3928 \\
   -1.6250 &  -1.1124 &  -0.8750 &  -0.0518 &  -1.1250  &  1.3624 &   2.1250  &  0.3018 \\
    1.3624 &  -1.1428 &  -0.1376 &   0.4608 &  -1.1982 &  -0.2286 &  -0.1124  &  0.5821 \\
    0.1250 &  -0.4053 &  -1.6250 &  -0.4053 &   2.1250  &  0.6553 &  -0.1250  &  0.6553 \\
    2.2160  &  0.6679 &  -1.4660 &  -3.1428  & -0.7589  &  0.6679  &  1.9231  & -0.5214 \\
   -0.1250  &  0.4482 &   0.1250 &  -0.9053 &  -1.6250 &   0.8018 &   0.1250  &  0.1553 \\
   -1.1124  &  0.4786 &  -2.6124 &  -0.8321 &  -1.5518 &   2.3928 &   2.3624  &  3.2892 \\
\end{array}\right)
$$

\end{document}